\documentclass[11pt,letterpaper]{amsart}
\usepackage{amsmath}
\usepackage{amssymb}

\usepackage{epsfig}
\usepackage[all, knot]{xy}
\xyoption{arc}

\newtheorem{thm}{Theorem}[section]
\newtheorem{df}[thm]{Definition}
\newtheorem{prop}[thm]{Proposition}
\newtheorem{cor}[thm]{Corollary}

\newtheorem{lem}[thm]{Lemma}
\newtheorem{ex}[thm]{Example}
\newtheorem{rem}[thm]{Remark}

\newtheorem*{thmA}{Theorem~A}
\newtheorem*{thmB}{Theorem~B}
\newtheorem*{thmC}{Theorem~C}
\newtheorem*{thmD}{Theorem~D}
\newtheorem*{thmE}{Theorem~E}

\newcommand{\Pic}{\operatorname{Pic}}

\newcommand{\ord}{\operatorname{ord}}
\newcommand{\Pre}{\operatorname{Preper}}

\newcommand{\pp}{\mathbb{P}}
\newcommand{\af}{\mathbb{A}}
\newcommand{\ox}{\mathcal{O}}
\newcommand{\Rat}{\operatorname{Rat}}

\newcommand{\afe}{\operatorname{AE}}

\begin{document}

\title[D-ratio and an upper bound for the height]
{Maximal ratio of coefficients of divisors \\
and an upper bound for height for rational maps}

\author[1]{Chong Gyu Lee}


\date{\today}


\address{Department of Mathematics, University of Illinois at Chicago, Chicago, IL, 60607, US}

\email{phiel@math.uic.edu}

\maketitle

\footnotetext[1]{This author was partially supported by the National Science Foundation and by the Korea Research Foundation funded by the Korea Government(MOEHRD). (KRF-2005-215-C00001) }

\begin{abstract}

    In this paper, we introduce the $D$-ratio of a rational map $f:\pp^n \dashrightarrow \pp^n$, defined over a number field $K$, whose indeterminacy locus is contained in a hyperplane $H$ on $\pp^n$. The $D$-ratio $r(f)$ provides useful height inequalities on $\pp^n(\overline{K}) \setminus H$: there is a constant $C$, depending only on $f$, such that
    \[
    \dfrac{r(f)}{\deg f} h\bigl( f(P) \bigr) + C > h(P)~\quad \text{for all}~P\in \pp^n(\overline{K}) \setminus H.
    \]
    If the indeterminacy loci of $f_1, f_2$ are disjoint subsets in $H$, then there is there is a constant $C'$, depending only on $f_1, f_2$, such that
    \[
        \dfrac{1}{\deg f_1} h\bigl( f_1(P) \bigr) + \dfrac{1}{\deg f_2} h\bigl( f_2(P) \bigr) +C'
        > \left( 1+ \min_{l=1,2} \left( \dfrac{1}{r(f_l)} \right) \right)h(P)
        \quad
        \text{for all }P \in \pp^n(\overline{K}) \setminus H.
    \]
    Also, we provide some dynamical applications of those height inequalities.
\end{abstract}

\section{Introduction}

    Let $f: \pp^n \dashrightarrow \pp^n$ be a rational map, defined over a number field $K$, and suppose its indeterminacy locus $I(f)$ is contained in a hyperplane $H$ on $\pp^n$. In this article, we introduce the $D$-ratio $r(f)$ associated to $f$ and a resolution of indeterminacy. We use $r(f)$ to provide a relation between $h(P)$ and $h\bigl(f(P) \bigr)$, where $h : \pp^n (\overline{K}) \rightarrow \mathbb{R}$ be the logarithmic absolute height function.
    \begin{thmA}
        Let $f :\pp^n \dashrightarrow \pp^n$ be a rational map, defined over a number field $K$, such that $I(f)$ is contained in a hyperplane $H$ and let $r(f)$ be a $D$-ratio of $f$ obtained from a resolution of indeterminacy. Then, there is a constant $C$, depending only on $f$, such that
        \[
        \dfrac{r(f)}{\deg f} h\bigl( f(P) \bigr) +C > h(P)
        \quad
        \text{for all }P\in \pp^n(\overline{K}) \setminus H.\]
    \end{thmA}

    \begin{thmB}
        Let $f_1, f_2:\pp^n \dashrightarrow \pp^n$ be rational maps, defined over a number field $K$, such that the indeterminacy loci of $f_1,f_2$ are disjoint subsets of a hyperplane $H$ and let $r(f_l)$ be a $D$-ratio of $f_l$. Then, there is a constant $C$, depending only on $f_1, f_2$, such that
        \[
        \dfrac{1}{\deg f_1} h\bigl( f_1(P) \bigr) + \dfrac{1}{\deg f_2} h\bigl( f_2(P) \bigr) +C
        > \left( 1+ \min_{l=1,2} \left( \dfrac{1}{r(f_l)} \right) \right)h(P)
        \quad
        \text{for all }P \in \pp^n(\overline{K}) \setminus H.\]
    \end{thmB}

     We have Northcott's theorem for endomorphisms: If $\phi:\pp^n \rightarrow \pp^n$ is an endomorphism, then there are two nonnegative constants $C_1, C_2$, depending on $\phi$, such that
        \begin{equation}\label{Northcott}\tag{A}
            \dfrac{1}{\deg \phi} h\bigl( \phi(P) \bigr) +C_1 > h(P)
            > \dfrac{1}{\deg \phi} h\bigl( \phi(P) \bigr) - C_2
        \quad \text{for all }P \in \pp^n(\overline{K}).
        \end{equation}
     Northcott's theorem is one of the essential theorems in arithmetic dynamics. For example, the Call-Silverman canonical height function \cite{CS} for an endomorphism on a projective space is well defined because of Northcott's theorem.

    In fact, Northcott's theorem only holds for endomorphisms:
    \begin{thmC}
        Let $f :\pp^n \dashrightarrow \pp^n$ be a rational map defined over a number field $K$.
        Suppose that $f$ satisfies the following inequality for some nonempty Zariski open set $U$ of $\pp^n$ and some constant $C$:
        \begin{equation*}
        \dfrac{1}{\deg f}h \bigl( f(P) \bigr) +C > h(P) \quad \text{for all }P\in U(\overline{K}).
        \end{equation*}
        Then, $f$ is an endomorphism.
    \end{thmC}
    \noindent Note that $h(P)$ is always bounded below by $\dfrac{1}{\deg f} h\bigl( f(P) \bigr) - C_2$ for any rational map $f$ of degree~$d$ so that we only care bout upper bound of $h(P)$. (For details, see \cite[Theorem~B.2.5]{SH}.)

     Thus, the canonical height functions and other arithmetic-dynamical consequences may fail to exist in the presence of indeterminacy.
     However, $D$-ratios provide weaker height inequalities as stated in Theorem~A and Theorem~B. Those inequalities allow us to study arithmetic-dynamical properties of rational maps. For example, the H\'{e}non map $g: K^n \rightarrow K^n$ is a polynomial map of special dynamical interest. It is not an endomorphism on a projective space and hence Northcott's theorem does not hold for $g$.
     Nevertheless, it is well known that the H\'{e}non map has good dynamical properties such as the boundedness of the set of periodic points. (For details, see \cite{D, M, S1}.) Especially, a regular H\'{e}non map $g$, together with the inverse $g^{-1}$, satisfies the hypothesis of Theorem~B and hence $g$ has a canonical height \cite{K, K2, Le}.

     The main idea of the $D$-ratio is to generalize the case of endomorphisms. The degree of an endomorphism $\phi :\pp^n \rightarrow \pp^n$ is the coefficient of $H$ in $\phi^*H=\deg \phi \cdot H$. In other words, we have
    \[
    \deg \phi = \sup \left\{ \delta ~\left|~ \dfrac{1}{\delta} \phi^*H -  H \text{ is effective in } \Pic(\pp^n)\otimes \mathbb{R} \right. \right\}.
    \]
    Then, the functorial property of the Weil height machine \cite[Theorem~B.3.2]{SH} gives the comparison of $h(P)$ and $h\bigl(\phi(P)\bigr)$:
    \[
    h_H \bigl( \phi^*(P) \bigr) = h_{\phi^*H}(P) + O(1) = \deg \phi \cdot h_H(P) + O(1).
    \]

    Let $f:\pp^n \dashrightarrow \pp^n$ be a rational map. Due to failure of the functoriality of the Weil height machine, we pass to a resolution of indeterminacy to work with a morphism: there exist a nonsingular projective variety $V$ and a birational morphism $\pi : V\rightarrow \pp^n$ such that $\widetilde{f} =f \circ \pi$ extends to a morphism.
    \[
        \xymatrix{
            V \ar[d]_{\pi} \ar[rd]^{\widetilde{f}}
            &  \\
            \pp^n \ar@{-->}[r]_{f} & \pp^n
         }
         \quad \quad
        \xymatrix{
            \pp^n \ar[d]_{id} \ar[rd]^{\phi} & \\
            \pp^n \ar[r]_{\phi} & \pp^n
        }
    \]
    (In the case of endomorphisms, we may think $\pi$ to be the identity map on $\pp^n$.)
    Now, we have two morphisms $\widetilde{f}, \pi :V \rightarrow \pp^n$. We can compare $\widetilde{f}^*H$ and $\pi^*H$ in $\Pic(V)$ and hence
    \[
    h_H\bigl(\widetilde{f}(P) \bigr) = h_{\widetilde{f}^*H}(P) +O(1)\quad \text{and} \quad h_H\bigl( \pi(P)\bigr) = h_{\pi^*H}(P) +O(1).
    \]
    Roughly, we define the $D$-ratio to be the constant $r(f)$ such that
    \begin{eqnarray*}\label{r(g)}
    \dfrac{\deg f}{r(f)} = \sup \left\{ \delta ~\left|~ \dfrac{1}{\delta}\widetilde{f}^*H - \pi^*H \text{ is }\af^n\text{-effective in } \Pic(V)\otimes \mathbb{R} \right. \right\}.
    \end{eqnarray*}
    Using properties of ``$\af^n$-effective'' divisors, we get Theorem~A and Theorem~B. Also, we conclude that
    \[
    r(f) = 1\quad \text{if and only if}\quad f~\text{is an endomorphism.}
    \]
    Note that the definition of the $D$-ratio depends on the choice of a resolution of indeterminacy. However, we show that it depends only on the ``strong factorization class'' of the resolution. (See Lemma~\ref{invariant}.) In particular, in dimension 2, the $D$-ratio depends only on $f$.

    We provide two applications of Theorem A and Theorem B in arithmetic dynamics. For convenience, define
    \[
    \Rat^n(H) := \left\{ f :\pp^n \dashrightarrow \pp^n ~\left| \right. I(f) \subset H \right\}.
    \]
    If $f\in \Rat^n(H)$ is a rational map such that $r(f) <\deg f$, then Theorem~A directly induces the following result:
    \begin{thmD}\label{preperiodic}
        Let $f:\af^n \rightarrow \af^n$ be a polynomial map, defined over a number field $K$, such that $r(f)<\deg f$. Then, the set of preperiodic points
        \[
        \Pre(f) := \left\{ P\in \af^n(\overline{K}) ~|~ f^l(P) = f^m(P) ~\text{for some}~l\neq m\right\}
        \]
        is a set of bounded height. Hence
        \[ \Pre(f) \cap \af^n(K') \]
        is a finite set for any number field $K'$.
    \end{thmD}

     The requirement $r(f)<\deg f$ in Theorem~D is sharp; there are rational maps such that $r(f) = \deg f$ and $\Pre(f)$ is not bounded. (See Example~\ref{RFD}.) Still, we can find some information for such rational map $f$ if it has a good counterpart.
     \begin{thmE}
        Let $S=\{f_1, f_2 \}$ be a pair of polynomial maps, defined over a number field $K$, such that their indeterminacy loci are disjoint subsets of a hyperplane $H$, let $f(f_l)$ be a $D$-ratio of $f_l$ and let $\Phi_S$ be the monoid of rational maps generated by $S$:
        \[
        \Phi_S := \{ f_{i_1} \circ \cdots \circ f_{i_m} ~|~ i_j = 1 ~\text{or}~2,~ m\geq 0\}.
        \]
        Define
        \[
        \delta_S : = \left(\dfrac{1}{1+1/r} \right) \left(\dfrac{1}{\deg f_1} + \dfrac{1}{\deg f_2} \right)
        \]
        where $r = \displaystyle  \max_{l=1,2} \bigl(r(f_l)\bigr)$.

        If $\delta_S <1$, then
        \[
        \Pre(\Phi_S):= \bigcap_{f \in \Phi_S} \Pre(f) \subset \af^n_{\overline{K}}
        \]
        is a set of bounded height.
    \end{thmE}

\par\noindent\emph{Acknowledgements}.\enspace
This paper is a part of my Ph.D. dissertation. I would like to thank my advisor Joseph H. Silverman for his overall advice. Also, thanks to Dan Abramovich for his helpful comments, especially for Section~2, and thanks to Laura DeMarco for useful discussions.

\section{Preliminaries : Blowup and its Picard Group}

      In this section, we check the basic theory of the resolution of indeterminacy. For details, I refer \cite{Cu,H} to the reader. We will let $H$ be a fixed hyperplane of $\pp^n$, let $\af^n = \pp^n \setminus H$ and let $f$ be an element of $\Rat^n(H)$ defined over a number field $K$ unless stated otherwise.

    \begin{thm}[Resolution of indeterminacy]
        Let $f : X \dashrightarrow Y$ be a rational map between proper varieties such that $X$ is nonsingular. Then there is a proper nonsingular variety $\widetilde{X}$ with a birational morphism $\pi: \widetilde{X} \rightarrow X$ such that $\phi = f \circ \pi : \widetilde{X} \rightarrow Y$ extends to a morphism:
        \[
                \xymatrix{
                \widetilde{X} \ar[d]_{\pi} \ar[rd]^{\phi} & \\
                X \ar@{-->}[r]_{f} & Y
                }
        \]
    \end{thm}

    Using Hironaka's Theorem (Theorem~\ref{Hironaka}), we will observe the relation between the resolution of indeterminacy and the indeterminacy locus of $f$.

    \begin{df}
        Let $\pi :\widetilde{X} \rightarrow X$ be a birational morphism. Then, we say that a closed subscheme $\mathfrak{I}$ of $X$ is \emph{the center scheme of $\pi$} if the ideal sheaf $\mathcal{S}$, corresponding to $\mathfrak{I}$, generates $\widetilde{X}$:
        \[
        \widetilde{X} = \operatorname{Proj} \left( \bigoplus_{d\geq 0}\mathcal{S}^d \right).
        \]
    \end{df}

    \begin{df}
        Let $\pi : \widetilde{X} \rightarrow X$ be a birational morphism. We say that $\pi$ is a \emph{monoidal transformation} if its center scheme is a smooth irreducible subvariety of $X$. We say that \emph{$\widetilde{X}$ is a successive blowup of $X$} if the corresponding birational map
        $\pi :\widetilde{X} \rightarrow X$ is a composition of monoidal transformations.
    \end{df}

    \begin{thm}[Hironaka]\label{Hironaka}
        Let $f :X \dashrightarrow Y$ be a rational map between proper varieties such that $X$ is nonsingular. Then, there is a finite sequence of proper varieties $X_0, \cdots, X_{r}$ such that
        \begin{enumerate}
        \item $X_0 = X$.
        \item $\rho_i : X_{i} \rightarrow X_{i-1}$ is a monoidal transformation.
        \item If $T_i$ is the center scheme of $\rho_i$, then $\rho_1 \circ \cdots \circ \rho_{i-1} (T_i) \subset I(f)$ on $X$.
        \item $f$  extends to a morphism $\widetilde{f}:X_{r} \rightarrow Y$.
        \item Consider the composition of all monoidal transformation $\rho : X_r \rightarrow X$. Then, the underlying subvariety of the center scheme $T$ of $\rho$, a subvariety made by the zero set of the ideal sheaf corresponding to $T$, is exactly $I(f)$.
        \end{enumerate}
    \end{thm}
    \begin{proof}
        See \cite[Question (E) and Main Theorem II]{Hi}.
    \end{proof}

    For notational convenience, we will define the following.
    \begin{df}
        Let $f :\pp^n \dashrightarrow \pp^n$ be a rational map.
          We say that a pair $(V,\pi)$ is \emph{a resolution of indeterminacy of $f$}
          when $V$ is a successive blowup of $\pp^n$ with a birational morphism $\pi :V \rightarrow \pp^n$ such that
        \[
        f \circ \pi : V \rightarrow \pp^n
        \]
        extends to a morphism. And we call the morphism $\phi := f \circ \pi$ a \emph{resolved morphism of $f$}.
    \end{df}

    In Section~3, we will find a basis of $\Pic(V)$ when $(V,\pi)$ is a resolution of indeterminacy of some rational map $f$. Especially, we need a basis consisting of irreducible divisors. However, pullbacks of irreducible divisors may not be irreducible because of the exceptional part. So, we define the proper transformation, which preserves the irreducibility.

    \begin{df}
        Let $\pi : \widetilde{X} \rightarrow X$ be a birational morphism with the center scheme $\mathfrak{I}$ and let $D$ be an irreducible divisor on $X$. We define \emph{the proper transformation of $D$ by $\pi$} to be
        \[
        \pi^{\#}D = \overline{\pi^{-1} (D \cap U)}
        \]
        where $U= X \setminus Z \left( \mathfrak{I}\right)$ and $Z \left( \mathfrak{I}\right)$ is the underlying subvariety made by the zero set of the ideal corresponding to $\mathfrak{I}$.
    \end{df}

    \begin{prop}\label{Pic}
        Let $V$ be a successive blowup of $\pp^n$ with a birational morphism $\pi:V \rightarrow \pp^n$: there are monoidal transformations $\pi_i : V_i \rightarrow V_{i-1}$ such that $V_r = V$ and $V_0 = \pp^n$.
        Let $H$ be a hyperplane on $\pp^n$, let $F_i$ be the exceptional divisor of the blowup $\pi_i : V_i
        \rightarrow V_{i-1}$, let $\rho_i = \pi_{i+1} \circ \cdots \circ \pi_{r}$ and let $E_i =\rho_i^\# F_i$. Then, $\Pic(V)$ is a free $\mathbb{Z}$-module with a basis
        \[
        \{H_V = \pi^\#H, E_1, \cdots, E_r \}.
        \]
    \end{prop}
    \begin{proof}
    \cite[Exer.II.7.9]{H} shows that
    \[
    \Pic(\widetilde{X}) \simeq \Pic(X) \oplus \mathbb{Z}
    \]
    if $\pi:\widetilde{X} \rightarrow X$ is a monoidal transformation. More precisely,
    \[
    \Pic(\widetilde{X}) = \{ \pi^\# D + nE ~|~ D \in \Pic(X)\}
    \]
    where $E$ is the exceptional divisor of $\pi$ on $\widetilde{X}$. Apply it to each $\pi_i$ and get the desired result.
    \end{proof}

\section{$\af^n$-effective divisor}

    We roughly describe the $D$-ratio of $f$ as follows: (Precise definition of the $D$-ratio will be given in Section~4.) $r(f)$ is the constant such that
    \begin{eqnarray*}\label{r(g)}
    \dfrac{\deg f}{r(f)} := \sup \left\{ \delta ~\left|~ \dfrac{1}{\delta}\widetilde{f}^*H - \pi^*H \text{ is }\af^n\text{-effective in } \Pic(V)\otimes \mathbb{R} \right. \right\}
    \end{eqnarray*}
    where $(V,\pi)$ is a resolution of indeterminacy, $\widetilde{f}$ is a resolved morphism and ``$\af^n$-effective divisor'' is the main topic of this section. Like the case of endomorphisms, we may use ``effective'' instead of new term ``$\af^n$-effective.'' However, it is hard to describe the effective cone of $V$ even though $\Pic(V)$ is a free $\mathbb{Z}$-module. Moreover, we cannot control the base locus of all effective divisors. So, we will take the $\af^n$-effective cone $\afe(V)$ such that 1)~$\afe(V)$ is a simple subset of the effective cone and 2)~all element of $\afe(V)$ has the base locus outside of $\af^n$. Remind that $H$ is a fixed hyperplane of $\pp^n$, $\af^n = \pp^n \setminus H$ and $f$ is an element of $\Rat^n(H)$ defined over a number field $K$.

    \begin{df}
        Let $V$ be a successive blowup of $\pp^n$ with a birational morphism $\pi:V \rightarrow \pp^n$ such that the underlying set of the center scheme of $\pi$ is a subset of $H$, let $H$ be a fixed hyperplane of $\pp^n$ and let
        \[
        \Pic_{\mathbb{R}}(V) = \mathbb{R} H_V \oplus \mathbb{R}E_1 \oplus \cdots \oplus  \mathbb{R}E_r
        \]
        with the basis $\{H_V , E_1, \cdots , E_s\}$ described in Proposition~\ref{Pic}. We define \emph{the $\af^n$-effective cone of $V$} to be
        \[
        \afe(V) := \mathbb{R}^{\geq 0} H_V \oplus \mathbb{R}^{\geq 0} E_1 \oplus \cdots \oplus  \mathbb{R}^{\geq 0} E_r
        \]
        where $\mathbb{R}^{\geq 0}$ is the set of nonnegative real numbers.
        We say that an element $D \in \Pic_{\mathbb{R}}(V)$ is \emph{$\af^n$-effective}
        if $D$ is contained in $\afe(V)$ and denote it by
        \[
        D \succ 0.
        \]
        Moreover, on $\Pic_{\mathbb{R}}(V)$, we write
        \[
        D_1 \succ D_2
        \]
        if $D_1 - D_2$ is $\af^n$-effective.
    \end{df}
    Recall that $H$ is a fixed hyperplane. So, we have a fixed basis of $\Pic(V)$. It implies the representation of an element in $\Pic(V)$ is unique and hence the $\af^n$-effectiveness is well defined.

    $\af^n$-effective divisors have some useful properties. To show them, we need the following lemma which is also important to define the $D$-ratio of a rational map later.

    \begin{lem}\label{m_{i,j}}
        Let $\pi: V\rightarrow \pp^n$ and $\rho: W\rightarrow V$ be compositions of monoidal transformations such that the underlying sets of the center schemes of $\pi$ and $\pi \circ \rho$ are subsets of $H$, let $\{H_V, E_1, \cdots , E_r\}$ $\{H_W, F_1, \cdots, F_s\}$ are bases of $\Pic(V)$ and $\Pic(W)$ respectively, described in Proposition~\ref{Pic} and let
        \[
        \rho^* H_V = \rho^{\#}H_V + \sum_{j=1}^s m_{0j}F_j \quad \text{and} \quad \rho^* E_i = \rho^{\#}E_i + \sum_{j=1}^s m_{ij}F_j.
        \]
        Then,
        \[
        m_{ij} \geq 0\quad \text{for all}~i,j.
        \]
        Furthermore,
        \[
        \displaystyle \sum_{i=0}^r m_{ij} >0\quad \text{for all }j=1, \cdots, s.
        \]
    \end{lem}
    \begin{proof}
        Fix $j\in \{ 1, \cdots, s\}$. Since the pullback of $E_i$ by $\rho$ is defined
        \[
        \rho^*E_i = \rho^{-1}E_i
        \]
        where $\rho^{-1}E_i$ is the scheme theoretic preimage. So, if $\rho(F_j) \subset E_i$, then $m_{ij} >0$. Otherwise, $m_{ij}=0$.
        Furthermore, since the underlying set $Z$ of the center of blowup of $W$ is in $H_V \cup \left(\displaystyle \bigcup_{i=1}^r E_i \right)$ by assumption, an irreducible subset $\rho(F_j)$ of $Z$ should be contained in one of irreducible components of $H_V \cup \left(\displaystyle \bigcup_{i=1}^r E_i \right)$, which is either $H_V$ or some $E_i$. Therefore, $m_{ij}>0$ for at least one $i$.
    \end{proof}

    \begin{prop}\label{af effec prop} Let $V$ be a successive blowup of $\pp^n$ with a birational morphism $\pi:V \rightarrow \pp^n$ and let $D, D_1, D_2, D_3 \in \Pic_{\mathbb{R}}(V)$.\\
        {\rm(1)  (Effectiveness)} If $D$ is $\af^n$-effective, then $D$ is effective.\\
        {\rm(2)  (Boundedness)} If $D$ is $\af^n$-effective, then
        $h_D(P)$ is bounded below on
        \[
        \pi^{-1} \bigl(\af^n\bigr) = V \setminus \left( H_V \cup \left(\displaystyle \bigcup_{i=1}^r E_i \right) \right).
        \] \\
        {\rm(3) (Transitivity)} If $D_1 \succ D_2$ and $D_2 \succ D_3$ , then $D_1 \succ D_3$\\
        {\rm (4) (Functoriality)} If $\rho:W \rightarrow V$ is a monoidal transformation and $D_1 \succ D_2$, then
        $\rho^*D_1 \succ \rho^* D_2.$
    \end{prop}
    \begin{proof}
    (1) It is obvious since $\afe(V)$ is a subset of the effective cone of $V$.\\

    (2) Since $D$ is $\af^n$-effective, it is effective. By the positivity of the Weil height machine \cite[Theorem~B.3.2.(e)]{SH}, we get
     that
    \[h_D(P) > \ox (1)
    \quad\text{for all}~P \in V \setminus |D|\]
     where $|D|$ is the base locus of $D$.
    The base locus of $D$, the intersection of all effective $(n-1)$-cycles linearly equivalent to $D$, is contained in any effective cycle linearly equivalent to $D$. By assumption, $D \sim p_0 H_V + \displaystyle \sum_{i=1}^r p_i E_i$ for some nonnegative integers $p_i$'s and hence
    $|D| \subset H_V \cup \left( \displaystyle\bigcup_{i=1}^r E_i \right)$.
    Therefore,
    \[
    V \setminus  \left( H_V \cup \left( \bigcup_{i=1}^r E_i \right) \right) \subset V \setminus |D|.
    \]
    Furthermore, by assumption $I(f)\in H$ and Theorem~\ref{Hironaka},
    \[
    \pi^{-1}\bigl(H\bigr) =  H_V \cup \left( \bigcup_{i=1}^r E_i \right)
    \]
    and hence
    \[
    \pi^{-1}\bigl(\pp^n \setminus H\bigr) =  V \setminus  \left( H_V \cup \left( \bigcup_{i=1}^r E_i \right) \right)
    \]

    (3) If $D_1 \succ D_2$ and $D_2 \succ D_3$, then $D_1-D_2$ and $D_2-D_3$ are in $\afe(V)$. Since $\afe(V)$ is closed under addition by definition, $D_1-D_3 = (D_1-D_2) + (D_2-D_3) \in \afe(V)$.\\

    (4) Let
    \[\Pic(V) = \mathbb{Z} H_V \oplus \mathbb{Z}E_1 \oplus \cdots \oplus  \mathbb{Z}E_r,\]
     let $W$ be a blowup of $V$
    with a monoidal transformation $\rho :W \rightarrow V$. Then, $\Pic(W)$ is still a free $\mathbb{Z}$-module:
    \[
    \Pic(W) = \mathbb{Z} H_V^\# \oplus \mathbb{Z} E_1^\# \oplus \cdots \oplus \mathbb{Z} E_r^\#
              \oplus \mathbb{Z} F
    \]
     where $H_V^\# = \rho^{\#}H_V$, $E_i^\#=\rho^{\#}E_i$ and
     $F$ is the exceptional divisor of $W$ over $V$. Moreover,
    \[
    \rho^*H_V = H_V^\# + m_{0}F \quad \text{and} \quad \rho^*E_i = E_i^\# + m_{i}F
    \]
    for some $m_{i}$, which are nonnegative integers by Lemma~\ref{m_{i,j}}.

    Therefore, for any $\af^n$-effective divisor $D=  p_0 H_V + \displaystyle \sum_{i=1}^r p_i E_i  \in \Pic_{\mathbb{R}}(V)$,
    \[
    \rho^*D = p_0 (\rho^*H_V) + \sum_{i=1}^r p_i (\rho^*E_i)
    =  p_0 H_V^\# + \sum_{i=1}^r p_i E_i^\# + \left( \sum_{i=0}^r p_i m_{i} \right) F
    \]
    is $\af^n$-effective on $W$ because $p_i$'s and $m_{i}$'s are nonnegative integers.
    \end{proof}

\section{Maximal ratio of coefficient of divisors}

    In this section, we introduce the main idea of this paper - the $D$-ratio.
    Since we fixed a hyperplane $H$ of $\pp^n$, we have a fixed basis of $\Pic(V)$ so that the representation of $D\in \Pic(V)$ is unique.
    Hence, the maximal ratio of coefficients of $H_V$ and $E_i$'s in $\phi^*H$ and $\pi^*H$ is well defined.
   \begin{df}
        Let $f \in \Rat^n(H)$, let $(V, \pi)$ be a resolution of indeterminacy of $f$ and let $\phi$ be the resolved morphism of $f$ on $V$:
            \[
                \xymatrix{
                V \ar[d]_{\pi} \ar[rd]^{\phi}\\
                \pp^n \ar@{-->}[r]_{f}  & \pp^n }
            \]
        Suppose that
            \[
            \pi^*H = a_0 H_V + \sum_{i=1}^r a_i E_i \quad \text{and} \quad \phi^*H = b_0 H_V + \sum_{i=1}^r b_i E_i.
            \]
        If $b_i$ are nonzero for all $i$ satisfying $a_i \neq 0$, we define \emph{the $D$-ratio} of $\phi$ to be
        \[r(\phi) := {\deg \phi} \cdot \max_i \left( \dfrac{a_i}{ b_i} \right).\]
        If there is an $i$ satisfying $a_i\neq 0$ and $b_i=0$, define
        \[r(\phi) := \infty.\]
    \end{df}

    \begin{rem}
        Let $f \in \Rat^n(H)$ and let $(V,\pi_V)$ be a resolution of indeterminacy of $f$. Then,
        \[
        {r(\phi)} : = \min  \left\{ C ~\left| \dfrac{C}{\deg f} \cdot \phi^*H - \pi^*H \succ 0 \right. \right\}.
        \]
    \end{rem}

   \begin{df}
        Let $f \in \Rat^n(H)$ be a rational map defined over a number field $K$. Then, we define \emph{the $D$-ratio of $f$ on $V$},
        \[
        r(f) = r\left( f,(V,\pi) \right) := r(\phi)
        \]
        where $(V,\pi)$ is a resolution of indeterminacy of $f$ described in Theorem~\ref{Hironaka} and $\phi$ is the resolved morphism of $f $ on $V$.
    \end{df}

    In the rest of this paper, we will use $r(f)$ instead of $r\left( f,(V,\pi) \right)$. The main application of the $D$-ratio is constructing height inequalities on $\af^n(\overline{K})$ and they do not care about the choice of $(V,\pi)$. Also, once we fix a resolution of indeterminacy $(V,\pi)$ of $f$, any successive blowup of $V$ will provides the same result.

    \begin{lem}\label{invariant}
        Let $(V,\pi_V)$ and $(W,\pi_W)$ be resolutions of indeterminacy of $f$ with resolved morphisms $\phi_V = f \circ \pi_V$ and $\phi_W = f \circ \pi_W$ respectively. Let $\tau = \pi_V^{-1} \circ \pi_W$.
            \[
                \xymatrix{
                W \ar[d]_{\pi_W}  \ar[rd]^{\phi_W} \ar@{-->}[rr]^{\tau} & &  V \ar[d]^{\pi_V} \ar[ld]_{\phi_V} \\
                \pp^n \ar@{-->}[r]_{f} & \pp^n   & \pp^n \ar@{-->}[l]^{f}
                }.
            \]
            Suppose that $\tau:W \dashrightarrow V$ allows strong factorization: there is a common blowup $U$ of $V$ and $W$ such that $\tau_V :U\rightarrow V$ and $\tau_W :U \rightarrow W$ are compositions of monoidal transformations.
            \[
                \xymatrix{
                 & U \ar[dd]|{\phi_U} \ar[rd]^{\tau_V} \ar[ld]_{\tau_W} &\\
                W \ar[d]_{\pi_W}  \ar[rd]^{\phi_W} & &  V \ar[d]^{\pi_V} \ar[ld]_{\phi_V} \\
                \pp^n \ar@{-->}[r]_{f} & \pp^n   & \pp^n \ar@{-->}[l]^{f}
                }.
            \]
             Then,
        \[
        r(\phi_V)= r(\phi_W).
        \]
        \end{lem}
        \begin{proof} Suppose
        \[
        \pi_V^*H = a_0 H_V + \sum_{i=1}^r a_i E_i\quad \text{and} \quad \phi_V^*H = b_0 H_V + \sum_{i=1}^r b_i E_i.
        \]

    First, consider the case that $W$ is a successive blowup of $V$. Suppose that $\rho :W \rightarrow V$ is a composition of monoidal transformations:
    \[
                \xymatrix{
                \pp^n \ar@{-->}[d]_{f} & V \ar[l]_{\pi_V} \ar[ld]|{\phi_V} & W \ar[l]_{\rho} \ar[lld]^{\phi_W}\\
                \pp^n &  &}.
     \]
    Since $\Pic(V) = \mathbb{Z} H_V \oplus \mathbb{Z}E_1 \oplus \cdots \oplus  \mathbb{Z}E_r$, we get
    \[
    \Pic(W) = \mathbb{Z} H_V^\# \oplus \mathbb{Z} E_1^\# \oplus \cdots \oplus \mathbb{Z} E_r^\#
              \oplus \mathbb{Z} F_1 \oplus \cdots \oplus \mathbb{Z} F_s
    \]
   where $H_V^\#= \rho^{\#}H_V$, $E_i^\#=\rho^{\#}E_i$ and $F_j$ are the irreducible components of the exceptional divisor of $W$ over $V$.
   Moreover, we may assume that
    \[
    \rho^*H_V = H_V^\# + \sum_{j=1}^s m_{0,j}F_j \quad \text{ and } \quad \rho^*E_i = E_i^\# + \sum_{j=1}^s m_{i,j}F_j
    \]
    for some integers $m_{i,j}$, which are nonnegative by Lemma~\ref{m_{i,j}}.
    By assumption, $\phi_W = \phi_V \circ \rho$ and hence
    \[
    \pi_W^*H = \rho^*\pi^*H = \rho^*\left( a_0 H_V + \sum_{i=1}^r a_i E_i \right) = a_0H_V^\# + \sum_{i=1}^r a_i E_i^\# + \sum_{j=1}^s \left(\sum_{i=0}^r a_im_{i,j} \right) F_j
    \]
    and
    \[
    \phi_W^*H = \rho^* \phi_V^*H = \rho^*\left(b_0 H_V + \sum_{i=1}^r b_i E_i \right) =b_0H_V^\# + \sum_{i=1}^r b_i E_i^\# + \sum_{j=1}^s \left(\sum_{i=0}^r b_im_{i,j} \right) F_j.
    \]
    If $b_i=0$ and $a_i \neq 0$ for some $i$, then $r(\phi_V)=r(\phi_W)=\infty$ by definition. So, we may assume $b_i>0$ for all $i$. Then, by definition of the $D$-ratio, we get an inequality
    \begin{equation}\label{aibi}\tag{B}
    r(\phi_V) = \deg \phi_V \cdot \max_i \left( \dfrac{a_i}{b_i}\right)\geq \deg \phi_V \cdot \dfrac{a_i}{b_i} \quad \text{for all}~i.
    \end{equation}
    Because Lemma~\ref{m_{i,j}} and the fact $b_i>0$, we have
    \[
    \sum_{i=0}^r b_i m_{i,j} \geq \sum_{i=0}^r m_{i,j} >0 \quad \text{for all }j,
    \]
    and hence all coefficient of $\phi_W^*H$ is positive. Thus, we have
    \[
    r(\phi_W) = \deg \phi_W \cdot \max \left(\max_i \left(\dfrac{a_i}{b_i}\right) , \max_j \left( \dfrac{\sum_{i=0}^r a_im_{i,j} }{\sum_{i=0}^r b_im_{i,j}} \right) \right).
    \]
    Moreover, due to (\ref{aibi}), we get
    \[
    \max_j \left( \dfrac{ \sum_{i=0}^r a_im_{i,j} }{ \sum_{i=0}^r b_im_{i,j} } \right)
    \leq \max_j \left( \dfrac{ \sum_{i=0}^r \frac{r(\phi_V)}{\deg \phi_V } b_im_{i,j} }{  \sum_{i=0}^r b_im_{i,j}  } \right) = \frac{r(\phi_V)}{\deg \phi_V } = \max_i \left(\dfrac{a_i}{b_i}\right).
    \]

    Finally, $\deg \phi_V = \deg \phi_W$ yields
    \[r(\phi_W) = \deg \phi_W \cdot \max \left(\max_i \left(\dfrac{a_i}{b_i}\right) , \max_j \left( \dfrac{\sum_{i=0}^r a_im_{i,j} }{\sum_{i=0}^r b_im_{i,j}} \right) \right)  = \deg \phi_V \cdot \max_i\left(\dfrac{a_i}{b_i}\right) = r(\phi_V).\]

    Now let $(V,\pi_V)$ and $(W,\pi_W)$ be resolutions of indeterminacy of $f$ allowing strong factorization: there is a common blowup $U$ of $V$ and $W$ such that $\tau_V :U\rightarrow V$ and $\tau_W :U \rightarrow W$ are compositions of monoidal transformation. Then, $(U, \pi_U:=\pi_V \circ \tau_V)$ is still a resolution of indeterminacy of $f$:
            \[
                \xymatrix{
                 & U \ar[dd]|{\phi_U} \ar[rd]^{\tau_V} \ar[ld]_{\tau_W} &\\
                W \ar[d]_{\pi_W}  \ar[rd]^{\phi_W} & &  V \ar[d]^{\pi_V} \ar[ld]_{\phi_V} \\
                \pp^n \ar@{-->}[r]_{f} & \pp^n   & \pp^n \ar@{-->}[l]^{f}
                }.
            \]
    Then, the previous result says
    \[
    r(\phi_V) = r(\phi_U) = r(\phi_W).
    \]
   \end{proof}

     \begin{prop}\label{D-ratio prop}
        Let $f,g\in \Rat^n(H)$ be rational maps defined over a number field $K$. Then,
        \begin{enumerate}
            \item $r(f)=1$ if $f$ is an endomorphism.
            \item $r(f) \geq 1 $ .
            \item There is a resolution of indeterminacy $(U,\pi_U)$ of $g\circ f$ such that the $D$-ratio of $g\circ f$ on $U$ satisfies the following inequality:
            \[
            \dfrac{r(f)}{\deg f} \cdot \dfrac{r(g)}{\deg g} \geq \dfrac{r(g\circ f)}{\deg (g\circ f) }.
            \]
            \item If $g$ is an endomorphism and $f$ is a rational map on $\pp^n$, then
                $   r(g\circ f)= r(f)$.
         \end{enumerate}
    \end{prop}
    \begin{proof} (1) When $f$ is an endomorphism, then $(\pp^n, id)$ is a resolution of indeterminacy of $f$. Thus,
    \[
    id^*H = H \quad \text{and} \quad f^*H = \deg f \cdot H
    \]
    and hence
    \[
    r(f) = \deg f \times \dfrac{1}{\deg f} = 1.
    \]
    If $(V,\pi)$ is an arbitrary resolution of indeterminacy of $f$, then $V$ is a successive blowup of $\pp^n$ so that
    \[
    r(f,(\pp^n,id)) = r(f,(V,\pi))
    \]
    because of Lemma~\ref{invariant}.

    (2) Let $(V,\pi)$ be a resolution of indeterminacy of $f$ with the resolved morphism $\phi=f \circ \pi$. We may assume that
    the underlying set of the center of blowup is $I(f)$ by Theorem~\ref{Hironaka}. Suppose that
    \[
    \pi^*H = a_0 H_V + \sum_{i=1}^r a_i E_i, \quad \phi^*H = b_0 H_V + \sum_{i=1}^r b_i E_i.
    \]

    We can easily check that $a_0=1$: because $\pi(E_i) \subset I(f)$ and $I(f)$ is a closed set of codimension at least $2$,
    ${\pi}_*E_i = 0$. Thus,
     \[
     {\pi}_*{\pi}^*H = {\pi_V}_*\left( a_0 H_V + \sum_{i=1}^r a_i E_i \right) = {\pi}_*a_0 H_V = a_0H.
     \]
     On the other hand, choose another hyperplane $H'$ which satisfies $I(f) \not \subset H'$. Then, since $\pi$ is one-to-one outside of $\pi^{-1}\left(I(f)\right)$, we have
    \[
    {\pi}_* \pi^*H = {\pi}_* \pi^*H'  = H' = H.
    \]
     Therefore, ${\pi}_*H_V = H$ and $a_0=1$.

    Now, let's figure $b_0$ out. We define the pull-back of $H$ by $\phi$ to be
   \begin{equation}\label{push}\tag{C}
   \phi^*H = \ord_{H_V}(u\circ \phi) \cdot H_V + \sum_{i=1}^r \ord_{E_i}(u\circ \phi) \cdot E_i
   \end{equation}
   where $u$ is a uniformizer at $H$. Apply $\pi_*$ on $(\ref{push})$ and get
   \[
   {\pi}_*\phi^*H = \ord_{H_V}(u \circ \phi) \cdot {\pi}_*H_V + \sum_{i=1}^r \ord_{E_i}(u\circ \phi) \cdot {\pi}_* E_i
   =\ord_{H_V}(u\circ \phi) \cdot H
   \]
   since ${\pi}_*H_V=H$ and ${\pi}_*E_i=0$.
   Furthermore, because $\phi = f \circ \pi $, $f =[x_0^d , f_1, \cdots, f_n]$ and $\pi$ is one-to-one on $H_V$, we get
   \[
   \ord_{H_V}(u \circ \phi) = \ord_{H}(u \circ \phi \circ \pi^{-1}) = \ord_{H}(u \circ f)= d.
   \]
   On the other hand, we have
   \[
   {\pi}_*\phi^*H = {\pi}_*\left( b_0H_V + \sum_{i=1}^r b_iE_i \right) = b_0 {\pi}_*H_V + \sum_{i=1}^r b_i {\pi}_*E_i = b_0H
   \]
   and hence $b_0=d$.

   Finally,
        \[
        r(f) = \deg f \cdot \max_i \left( \dfrac{a_i}{b_i} \right) \geq \deg f \cdot  \dfrac{a_0}{b_0} = \deg f \cdot \dfrac{1}{\deg f} = 1.
        \]

    (3)   If $r(f)$ or $r(g)= \infty$, then it is clear. So, we may assume that $r(f),r(g)$ are finite.
        Let $(V,\pi_V)$ and $(W,\pi_W)$ be resolutions of indeterminacy of $f$ and $g$ obtained by Theorem~\ref{Hironaka} respectively and
        $\phi = f \circ \pi_V, \psi = g \circ \pi_W$ are resolved morphisms defining the $D$-ratio of $f$ and $g$ respectively. Consider a rational map $\phi' = \pi_W^{-1} \circ \phi: V \dashrightarrow W$ and find a resolution of indeterminacy $(U, \rho)$ of $\phi'$ over $V$:
            \[
                \xymatrix{
                U \ar@/_2pc/[dd]_{\pi_U} \ar[d]_{\rho} \ar[rd]^{\widetilde{\phi}} & &\\
                V \ar[d]_{\pi_V} \ar@{-->}[r]^{\phi'} \ar[rd]^{\phi} &  W \ar[d]^{\pi_W} \ar[rd]^{\psi} &    \\
                \pp^n \ar@{-->}[r]_{f}  & \pp^n \ar@{-->}[r]_{g}  & \pp^n
                }
            \]
        Note that $U$ is a successive blowup of $V$ and hence $r(\phi) = r(\rho\circ \phi)$.

        Let
        \[
        \alpha = \dfrac{r(\phi)}{\deg \phi} \quad \text{and} \quad \beta = \dfrac{r(\psi)}{\deg \psi}.
        \]
        By the definition of $r(\phi)$ and $r(\psi)$, we have
        \[
        \pi_V^*H \prec \alpha \cdot\phi^*H \quad \text{and} \quad \pi_W^*H \prec \beta \cdot\psi^*H.
        \]
        By the functoriality of the $\af^n$-effectiveness, we get
        \begin{equation}\label{afine}\tag{D}
        \rho^*\pi_V^*H \prec \alpha \cdot\rho^* \phi^*H\quad \text{and} \quad\widetilde{ \phi}^*\pi_W^*H \prec \beta \cdot\widetilde{ \phi}^* \psi^*H.
        \end{equation}
        Since the diagram commutes: $ \phi\circ \rho =  \pi_W \circ\widetilde{ \phi}$, we get
        \[
        \rho^* \phi^*H =  \widetilde{ \phi}^*\pi_W^* H.
        \]
        So, we can connect inequalities in (\ref{afine}):
        \[
        \pi_U^*H = \rho^*\pi_V^*H  \prec \alpha \cdot \rho^*\phi^*H = \alpha \cdot\widetilde{\phi}^*\pi_W^* H  \prec \alpha \beta \cdot \widetilde{ \phi}^* \psi^*H.
        \]

        Therefore, $\alpha \beta$ is a constant satisfying
        \[\alpha \beta \cdot (\psi \circ \widetilde{ \phi} )^*H - \pi_U^* H =  \alpha \beta \cdot \widetilde{ \phi}^* \psi^*H -  \rho^* \pi_V^* H \succ 0,\]
        where $\psi \circ \widetilde{ \phi}$ is a resolved morphism of $g\circ f$.
        It follows that $\alpha\beta \geq \dfrac{r(\psi \circ \widetilde{ \phi})}{ \deg \psi \circ \widetilde{ \phi}}$. \\

    (4) Let $(V,\pi)$ be a resolution of indeterminacy of $f$ and suppose that
    \[
    \pi^*H = a_0 H_V + \sum_{i=0}^r a_i E_i, \quad \phi^*H = b_0 H_V + \sum_{i=0}^r b_i E_i.
    \]
    Consider the following diagram:
     \[
                \xymatrix{
                V \ar[d]_{\pi} \ar[rd]^{\phi}\\
                \pp^n \ar@{-->}[r]_{f} \ar[d]_{id_{\pp^n}}  & \pp^n  \ar[d]_{id_{\pp^n}} \ar[rd]^{g} \\
                \pp^n \ar@{-->}[r]_{f}  & \pp^n   \ar[r]_{g} & \pp^n
                }.
      \]
    Thus, we get
    \[
    \pi^* id_{\pp^n}^* H = \left( a_0 H_V + \sum_{i=1}^r a_i E_i \right) , \quad  \phi^*g^*H = \deg g \left( b_0 H_V + \sum_{i=1}^r b_i E_i\right)
    \]
    and hence
    \[
    \dfrac{r(g\circ f)}{\deg (g \circ f)} = \max_i \left(\dfrac{a_i }{\deg g \cdot b_i } \right) =  \dfrac{1}{\deg g} \max_i \left(\dfrac{ a_i }{b_i } \right)
    = \dfrac{r(f)}{\deg f \deg g}.
    \]
    Furthermore, $\deg (g \circ f) = \deg f \cdot \deg g$ since $g$ is an endomorphism. Therefore, we have the desired result:
    \[
    r(g\circ f) =r(f).
    \]
    \end{proof}

\section{Upper bounds for height for rational map}

    In this section, we prove Theorem~A and apply it to arithmetic dynamics. We start with Theorem~C, which says that we can only expect a weaker height inequality than Northcott's theorem.
    \begin{thmC}
        Let $f \in \Rat^n(H)$ be a rational map defined over a number field $K$. Suppose that $f$ satisfies the following inequality
        for some nonempty Zariski open set $U$ of $\pp^n$ and some constant $C$:
        \begin{equation}\label{ineq}\tag{E}
        \dfrac{1}{\deg f}h \bigl( f(P) \bigr) +C > h(P) \quad \text{for all }P\in U(\overline{K})\setminus I(f).
        \end{equation}
        Then, $f$ is an endomorphism.
    \end{thmC}
    \begin{proof}
    Suppose that there is a point $Q\in I(f)$. Without loss of generality, we may assume that $Q = [0,0, \cdots, 0, 1]$. Let
    \[
    f(X) = [f_0(X), f_1(X), \cdots f_n(X)]
    \]
    where $X=[X_0, \cdots, X_n]$, $d=\deg f$ and $f_i$ are homogeneous polynomials of degree $d$.
    Then, we can claim that
    \[
    \deg _{X_n} f_i < d
    \]
    for all $i=0, \cdots , n$: if there is an $j$ such that $\deg_{X_n}f_j = d$, then $f_j(Q) \neq 0$ and hence $Q$ cannot be an indeterminacy point. Choose a point $\alpha = [\alpha_0, \cdots, \alpha_{n-1}] \in \pp^{n-1}$ and
    a projective line
    \[
    L_{\alpha} := \left\{ [x_0, \cdots, x_n] ~|~ \alpha_i x_j = \alpha_j x_i \text{ for all } 0 \leq i,j \leq n-1 \right\} \not \subset H.
    \]
    Precisely, $L$ is the image of the closed embedding of $\pp^1$:
    \[
    \iota_{{\alpha}}: \pp^1\rightarrow L \subset \pp^n, \quad [Y_0, Y_1] \mapsto [\alpha_0 Y_0, \cdots, \alpha_{n-1} Y_0, Y_1].
    \]
    Since $U$ is dense, so is $U' = U \setminus I(f)$. Therefore, there is an $\alpha$ such that $L_\alpha \cap U'$ is dense in $L_\alpha$.
    Moreover, $f$ is defined on $L_\alpha \cap U' = \{ [\alpha_0,\cdots, \alpha_{n-1}, x_n]\}$ where $x_n = Y_1 / Y_0$ and hence
    $f_i[\alpha_0,\cdots, \alpha_{n-1}, x_n]$ is one-variable polynomial of degree at most $\deg _{X_n} f_i$ for all $i=0, \cdots, n$. Thus,
    $f|_L$ is a morphism of degree $d' \leq \max_{i} \deg _{X_n} f_i <d$ on $L$ and hence we have the following inequality:
    \[
    h(P) > \dfrac{1}{d'} h \bigl( f(P) \bigr) - C'
    \quad
    \text{for all }P\in L(\overline{K}) \cap U',\]
    which contradicts (\ref{ineq}).
    \end{proof}

    \begin{thmA}\label{main1}
    Let $f\in \Rat^n(H)$ be a rational map defined over a number field $K$ and let $r(f)$ be a $D$-ratio of $f$.
    Then, there is a constant $C$, depending only on $f$, such that
    \[
    \dfrac{ r(f) }{\deg f} h(f(P)) +C >  h(P)
    \quad
    \text{for all }P \in \af^n(\overline{K}).\]
    \end{thmA}
    \begin{proof}
    Let $(V,\pi)$ be a resolution of indeterminacy of $f$ with the resolved morphism $\phi = f \circ \pi_V$ and $r(f)$ be the $D$-ratio of $f$ on $V$. Suppose that
    \[
    \pi^*H = a_0 H_V +\sum_{i=1}^r a_i E_i \quad \text{and} \quad \phi^* H= b_0 H_V +\sum_{i=1}^r b_i E_i
    \]
    where $\{H_V, E_1, \cdots, E_r\}$ is the basis described in Proposition~\ref{Pic}.

    Let
    \[
    E := \dfrac{r(f)}{\deg f} \phi^*H - \pi^*H.
    \]
    By the definition of the $D$-ratio, $E$ is $\af^n$-effective. So, the height function corresponding to $E$,
    \[
    h_{E} = \dfrac{r(f)}{\deg f}h_{ \phi^*H} - h_{\pi^*H}
    \]
    is bounded below on $\pi^{-1}\bigl(\af^n\bigr)$ by Proposition~\ref{af effec prop}~(1).
    Hence we have the following inequality:
    \[
    \dfrac{r(f)}{\deg f} h\bigl( \phi(Q) \bigr) +C >  h\bigl( \pi_V(Q) \bigr) \quad \text{for all}~Q\in \pi^{-1}\bigl(\af^n\bigr).\]

    Let $P \in \af^n(\overline{K})$ be an arbitrary point and take $Q_P=\pi^{-1}(P)$. Then $\phi(Q_P)=f(P)$ and $Q_P$ satisfies the above inequality. Therefore, we get the desired result:
    \[
    \dfrac{r(f)}{\deg f} h\bigl( f(P) \bigr) +C >  h(P) \quad \text{for all }P \in \af^n(\overline{K}).
    \]
    \end{proof}

    \begin{cor}
        Let $f \in \Rat^n(H)$ be a rational map defined over a number field $K$. Then
        \[
        r(f) = 1 \quad \text{if and only if} \quad f \text{ is an endomorphism.}
        \]
    \end{cor}
    \begin{proof}
    One direction is already done; see Proposition~\ref{D-ratio prop} (1). For the other direction, suppose that $r(f)=1$. Then,
    by Theorem~A, we get the inequality $(E)$. So, $f$ is an endomorphism by Theorem~C.
    \end{proof}

    We apply Theorem~A to study dynamics of polynomial map. We can consider that a polynomial map $f:\af^n \rightarrow \af^n$ is an element $f \in\Rat^n(H)$ such that $f(\af^n) \subset \af^n$. Thus, we can define $r(f)$ and apply Theorem~A at all forward image $f^m(P)\in \af^n$.
    \begin{thmD}\label{preperiodic}
        Let $f:\af^n \rightarrow \af^n$ be a polynomial map, defined over a number field $K$, such that $r(f)<\deg f$. Then,
        \[
        \Pre(f) = \left\{ P\in \af^n(\overline{K}) ~|~ f^l(P) = f^m(P) ~\text{for some}~l,m\right\}
        \]
        is a set of bounded height and hence
        \[ \Pre(f) \cap \af^n(K') \]
        is finite for any number field $K'$.
    \end{thmD}
    \begin{proof}
    Let $u= \dfrac{r(f)}{\deg f} <1$. Then, by Theorem~A, we have
    \begin{equation}\label{basic}\tag{F}
    u\cdot  h\bigl( f(P) \bigr) >  h(P) - C \quad \text{for all }P \in \af^n(\overline{K}).
    \end{equation}
     Then, the iteration of (\ref{basic}) provides
    \begin{eqnarray*}
    u^{l} \cdot h\bigl( f^{l}(P) \bigr)
    &>& u^{l} \cdot h\bigl( f(f^{l-1}(P)) \bigr)\\
    &>& u^{l-1} \left[ h\bigl( f^{l-1}(P) \bigr) - C \right] \\
    &>& u^{l-2} \left[ h\bigl( f^{l-2}(P) \bigr) -C \right] - u^{l-1} C\\
    &>& \hspace{2cm} \vdots\\
    &>& h(P) - \left\{ 1+ u  + \cdots + u^{l-1}  \right\} C.\\
    \end{eqnarray*}
    Hence, we have
    \[
    \lim_{l \rightarrow \infty} u^{l} \cdot  h\bigl( f^{l}(P) \bigr)
    > h(P) - \dfrac{C}{1-u}.
    \]
    If $P$ is a preperiodic point of $f$, then the left hand side goes to zero so that
    \[\dfrac{C}{1-u}> h(P)\quad
    \text{for all }P \in \af^n(\overline{K}).\]
    \end{proof}

    \begin{ex}\label{ex1}
        Let
        \[
        F(x,y) = (x^3+y,x+y^2).
        \]
         Then, after three blowing-ups along points (see Figure~1), we get a resolution of indeterminacy of $F$. And, we have
         \[
         \pi^*H = H_V + E_1 + 2E_2 + 3E_3, \quad \phi^*H = 3H_V + 2E_1 + 4E_2 + 6E_3.
         \]
         Thus, $r(F) = 3/2<3$ and hence it has finitely many preperiodic points on $\af^n(K)$.
    \end{ex}
         \begin{figure}[h]
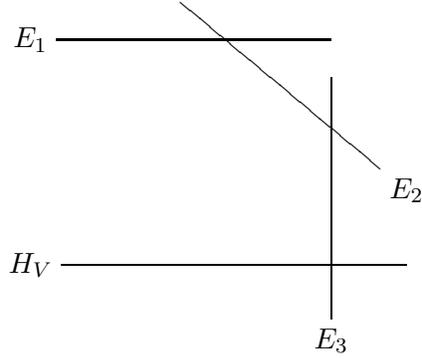

         \[
         \xy <1cm,0cm>:
        (0,0) *+{H_V}; (5,0) **@{-},
        (2,3.5) ; (5,1)*+{E_2} **@{-},
        (0,3)*+{E_{1}} ; (4,3) **@{-},
        (4,2.5) ; (4,-1)*+{E_{3}} **@{-},
        \endxy
        \]
        \caption{Resolution of indeterminacy of $F[X,Y,Z] = [X^3+YZ^2, XZ^2+Y^2Z,Z^3]$}
        \end{figure}

    \begin{ex}\label{RFD}
        The condition $r(f) < \deg f$ in Theorem~D is sharp: let
        \[
        F_0(x,y) = (x,y^2).
        \]
         Then, after two blowing-ups along points (see Figure~2), we get a resolution of indeterminacy of $F_0$. And, we have
         \[
         \pi^*H = H_V + E_1 + 2E_2, \quad \phi^*H = 2H_V + E_1 + 2E_2.
         \]
         Thus, $r(F_0) = 2 = \deg F_0$. And, it has infinitely many integral fixed points $(n,0)$. Thus, $\Pre(F_0)$ is not bounded.
    \end{ex}
         \begin{figure}[h]
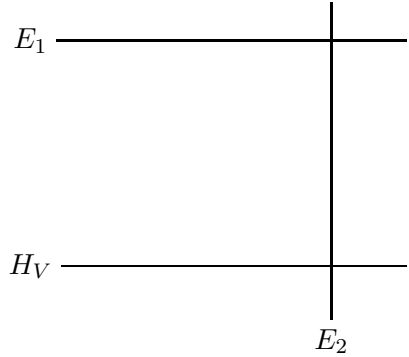

         \[
         \xy <1cm,0cm>:
        (0,0) *+{H_V}; (5,0) **@{-},
        (0,3)*+{E_{1}} ; (5,3) **@{-},
        (4,3.5) ; (4,-1)*+{E_{2}} **@{-},
        \endxy
        \]
        \caption{Resolution of indeterminacy of $F_0[X,Y,Z] = [XZ, Y^2,Z^2]$}
        \end{figure}

    \begin{cor}
        Let $f\in \Rat^n(H)$ be a rational map, defined over a number field $K$. If there is some number $N$ satisfying $r(f^N) < \deg (f^N)$, then
        \[\Pre(f) = \left\{ P\in \af^n ~|~ f^l(P) = f^m(P) ~\text{for some}~l\neq m\right\}\]
        is a set of bounded height.
    \end{cor}
    \begin{proof}
    It is enough to show that
    \[\Pre(f) = \Pre(f^N).\]
    One direction is clear;
    \begin{eqnarray*}
    P \in\Pre(f) &\Rightarrow& \ox_f(P) \text{ is finite} \\
                 &\Rightarrow&  \ox_{f^N} (P) \text{ is finite since}~ \ox_{f^n} (P) \subset \ox_{f} (P)\\
                 &\Rightarrow& p \in \Pre(f^N).
    \end{eqnarray*}
    So, we only have to show the other direction.
    Suppose that $P$ is a preperiodic point of $f$. Then,
    \[f^l(P) = f^m(P)\]
    holds for some natural numbers $l< m$. Let $g = f^{m-l}$. Then, we have
    \[g\bigl( f^l(P) \bigr) = f^l(P),\]
    which means $f^l(P)$ is a fixed point of $g$.
    So, we get
    \[
    g^N\bigl( f^l(P) \bigr) = f^l(P)
    \]
    and
    \[
    f^{kN-l}\bigl( g^N\bigl( f^l(P) \bigr) \bigr) = f^{kN-l} \bigl( f^l(P) \bigr) = f^{kN}(P)
    \]
    for all positive integers $k,N$. Since $g = f^{m-l}$, we get
    \[ f^{N(k+m-l)}(P) =  f^{kN -l +N(m-l) +l}(P) = f^{kN-l} (g^N(f^l(P)))  = (f^N)^k(P),\]
    which means that $P$ is a preperiodic point of $f^N$.
    \end{proof}

    \begin{ex}
        Let
        \[
        G(x,y) = (y,x^2+y).
        \]
         Then, after two blowing-ups (see Figure~3), the indeterminacy of $G$ is resolved.
         And, we have
         \[
         \pi^*H = H_V + E_1 + 2E_2 \quad \text{and} \quad \psi^*H = 2H_v + E_1 + 2E_2
         \]
         so that $r(G) = 2$.
         But $G^2(x,y) = (x^2+y,x^2+y^2+y)$ which extends to a morphism. Thus, $r(G^2)=1<3$. Therefore, it has finitely many preperiodic points on $\af^n(K)$.
    \end{ex}
        \begin{figure}[h]
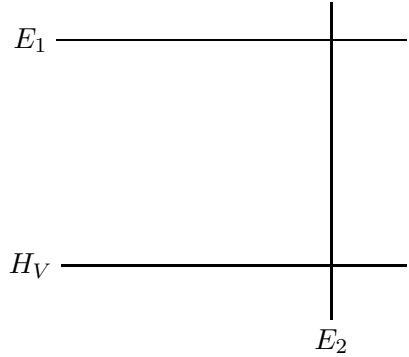

         \[
         \xy <1cm,0cm>:
        (0,0) *+{H_V}; (5,0) **@{-},
        (0,3)*+{E_{1}} ; (5,3) **@{-},
        (4,3.5) ; (4,-1)*+{E_{2}} **@{-},
        \endxy
        \]
        \caption{Resolution of indeterminacy of $G[X,Y,Z] = [YZ, X^2+YZ,Z^2]$}
        \end{figure}

    \begin{ex}
        Let \[f_{Na}(x,y,z) = (x+(x^2-yz)z, y+2(x^2-yz)x+(x^2-yz)^2z,z)\] be the Nagata map. Then, any points on the quadratic curve $x^2 = yz$ is a fixed point of $f_{Na}$. So, for any $N$, $r(f_{Na}^N) \geq \deg f_{Na}^N$. Furthermore, the $N$-th iteration of $f_{Na}$ is still a polynomial map of degree $5$;
        \[f_{Na}^N(x,y,z) = (x+N(x^2-yz)z, y+2N(x^2-yz)x+N^2(x^2-yz)^2z,z).\]
        Therefore, for any resolution of indeterminacy $(V_m,\pi_m)$ of any iteration of $f_{Na}^N$, the resolved morphism $\phi_N$ satisfies
        \[r(\phi_N) \geq 5 .\]
    \end{ex}

    \begin{rem}
        Theorem~D only considers preperiodic points of $f$ in $\af^n(\overline{K})$. Even though a rational map $f$ satisfying $r(f)<\deg f$ could have infinitely many periodic points on the hyperplane $H$.
    \end{rem}
    \begin{ex}
        Consider a rational map on Example~\ref{ex1}:
        \[
        F[x,y,z] = [x^3+ yz^2,xz^2+y^2z, z^3].
        \]
        It has infinitely many preperiodic points $P=[a,b,0]$ with $a \neq 0.$
    \end{ex}

\section{Jointly Regular Pairs}

     In this section, we will prove Theorem~B, which is a generalization of results of Silverman \cite{S4}, Kawaguchi \cite{K} and the author \cite{Le} for jointly regular pairs: we say that $S$ is a jointly regular pair if $S=\{ f_1, f_2 \}$ is a set consisting of two rational maps whose their indeterminacy loci are disjoint.

     Silverman \cite{S4} proved a weaker result for a jointly regular family of polynomial maps: Let $\{f_1, \cdots, f_k\}$
     be a family of polynomial maps, defined over a number field $K$. Suppose that the intersection of indeterminacy loci of $f_l$'s is empty. Then there is a constant $C$ satisfying
        \[
        \sum_{l=1}^k \dfrac{1}{\deg f_l} h\bigl( f_l(P) \bigr) +C > h(P)
        \quad
        \text{for all }P \in \pp^n(\overline{K}) \setminus H.\]
    Recently, Kawaguchi \cite{K2} and the author \cite{Le} independently proved Theorem~B for regular polynomial automorphisms: we say that \emph{$f$ is a regular polynomial automorphism} if $f:\af^n \rightarrow \af^n$ has the inverse map $f^{-1}:\af^n \rightarrow \af^n$ and $I(f)\cap I(f^{-1}) = \emptyset$. Then there is a constant $C$ satisfying
        \[
        \dfrac{1}{\deg f} h\bigl( f(P) \bigr) +\dfrac{1}{\deg f^{-1}} h\bigl( f^{-1}(P) \bigr) +C> \left(1+ \dfrac{1}{\deg f \deg f^{-1}} \right)h(P)
        \quad
        \text{for all }P\in \af^n(\overline{K}).\]
    Note that $r(f) = \deg f \cdot \deg f^{-1}$ if $f$ is a regular polynomial automorphism. (For details, see \cite{Le}.)

    \begin{thmB}
        Let $\{ f_1, f_2 \} \subset \Rat^n(H)$ be a jointly regular pair of rational maps defined over a number field $K$ and let $r(f_l)$ be the $D$-ratio of $f_l$. Then, there is a constant $C$ satisfying
        \[
        \dfrac{1}{\deg f_1} h\bigl( f_1(P) \bigr) + \dfrac{1}{\deg f_2} h\bigl( f_2(P) \bigr) +C
        > \left( 1+ \min_{l=1,2} \left( \dfrac{1}{r(f_l)} \right) \right)h(P)
        \quad
        \text{for all }P \in \pp^n(\overline{K}) \setminus H.\]
    \end{thmB}

    \begin{proof}
        Let $(V_l, \pi_l)$ be resolutions of indeterminacy of $f_l$ obtained by Theorem~\ref{Hironaka} and let
        $\phi_l = f_l \circ \pi_l$ be the resolved morphisms of $f_l$ respectively. Then, the underlying set of the center scheme of $V_l$ are exactly $I(f_l)$. So, we may assume that the centers of blowups $V_1,V_2$ are disjoint. Suppose that
        \[
        \Pic(V_l) = \mathbb{Z} \pi_l^\#H \oplus \mathbb{Z}E_{l1}\oplus \cdots\oplus \mathbb{Z}E_{lr_l}.
        \]
        and
        \[ \pi_l^*H = \pi_l^\# H + \sum_{i=1}^{r_l} a_{li} E_{li},
        \quad  {\phi_l}^*H = d_{l} \cdot \pi_l^\#H + \sum_{i=1}^{r_l} b_{li} E_{li}.\]
        where $d_l = \deg f_l$. Note that $b_{l0}=d_l$ from the proof of Proposition~\ref{D-ratio prop}~$(1)$.

        Now, consider a blowup $U$ of $\pp^n$ along union of centers of $V_1$ and $V_2$. Then, we have the following diagram:
        
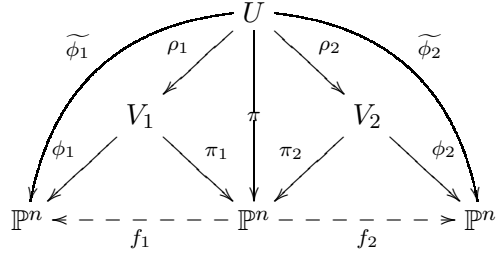
\[
            \xymatrix{
                & &  U \ar[dd]|{\pi}
                \ar[ld]_{\rho_1} \ar[rd]^{\rho_2} \ar@/_2pc/[lldd]_{\widetilde{\phi_1}} \ar@/^2pc/[rrdd]^{\widetilde{\phi_2}}& & \\
                & V_1 \ar[ld]_{\phi_1} \ar[rd]^{\pi_1} & & V_2 \ar[rd]^{\phi_2} \ar[ld]_{\pi_2} & \\
                \pp^n & & \pp^n \ar@{-->}[ll]^{f_1} \ar@{-->}[rr]_{f_2} & & \pp^n
            }
        \]
        Because $I(f_1)$ and $I(f_2)$ are disjoint, $U$ is still a successive blowup of $\pp^n$. (For details, See Figure 4.)
                \begin{figure}[b]
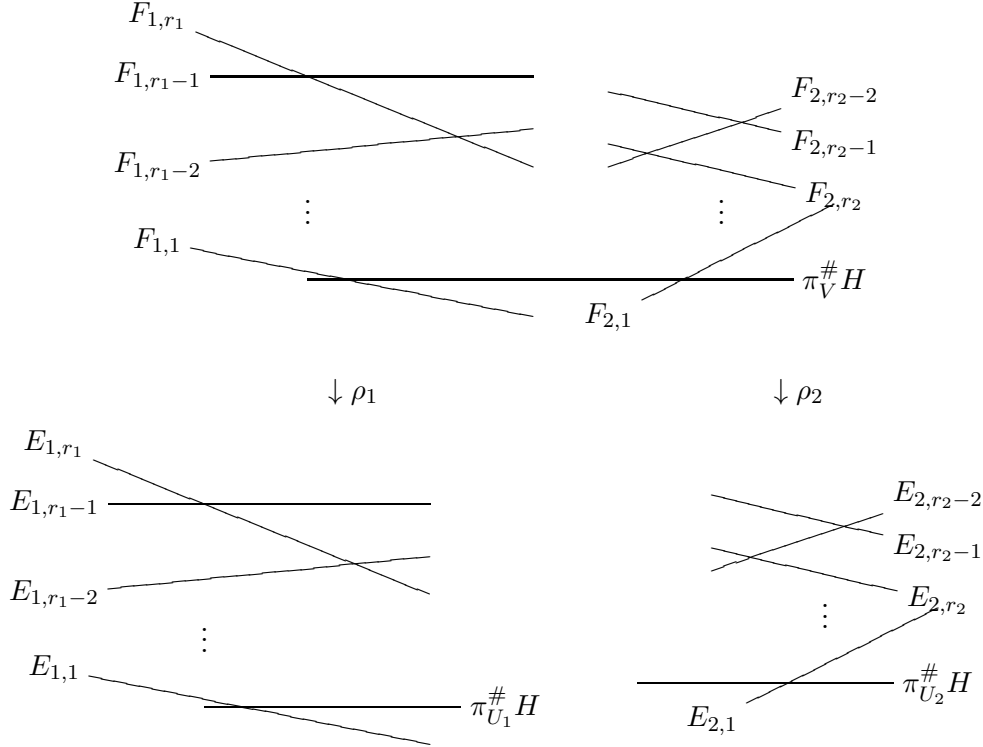

        \[
         \xy <1cm,0cm>:
        (0,3) *+{F_{1,1}}; (5,2) **@{-},
        (2,2.5) ; (9,2.5)*+{\pi_{V}^\#H} **@{-},
        (6,2)*+{F_{2,1}} ; (9,3.5)  **@{-},
        (7.5,3.5) *+{\vdots},
        (6,4) ; (9,5) *+{F_{2,r_2-2}} **@{-},
        (6,5) ; (9,4.3) *+{F_{2,r_2-1}} **@{-},
        (6,4.3) ; (9,3.6) *+{F_{2,r_2}} **@{-},
        (2,3.5) *+{\vdots},
        (0,5.2) *+{F_{1,r_1-1}}; (5,5.2) **@{-},
        (0,6) *+{F_{1,r_1}}; (5,4) **@{-},
        (0,4) *+{F_{1,r_1-2}}; (5,4.5) **@{-},
        \endxy
        \]
        \hspace*{2cm} $\downarrow\rho_1$ \hspace*{5cm} $\downarrow\rho_2$
        \[
           \xy <1cm,0cm>:
        (0,3) *+{E_{1,1}}; (5,2) **@{-},
        (2,2.5) ; (6,2.5)*+{\pi_{U_1}^\#H} **@{-},
        (2,3.5) *+{\vdots},
        (0,5.2) *+{E_{1,r_1-1}}; (5,5.2) **@{-},
        (0,6) *+{E_{1,r_1}}; (5,4) **@{-},
        (0,4) *+{E_{1,r_1-2}}; (5,4.5) **@{-},
        \endxy
        \quad \quad \quad
        \xy <1cm,0cm>:
        (5,2.5) ; (9,2.5)*+{\pi_{U_2}^\#H} **@{-},
        (6,2)*+{E_{2,1}} ; (9,3.5)  **@{-},
        (7.5,3.5) *+{\vdots},
        (6,4) ; (9,5) *+{E_{2,r_2-2}} **@{-},
        (6,5) ; (9,4.3) *+{E_{2,r_2-1}} **@{-},
        (6,4.3) ; (9,3.6) *+{E_{2,r_2}} **@{-},
        \endxy
        \]
        \caption{Exceptional divisors of $V_1$,$V_2$ and $U$}
        \end{figure}
        Moreover, we have
        \[
        \rho_1^* E_{1j} = \rho_1^\# E_{1j}\quad \text{and} \quad \rho_2^* E_{2j} = \rho_2^\# E_{2j}.
        \]
        Thus, let $F_{lj} = \rho_l^\# E_{lj}$ and get the following description of $\Pic(U)$:
        \[\Pic (U) = \mathbb{Z}\pi^\#H\oplus \mathbb{Z}F_{11}\oplus \cdots \oplus\mathbb{Z}F_{1r_1}\oplus
        \mathbb{Z}F_{21}\oplus \cdots \oplus\mathbb{Z}F_{2r_2}\]
        \[\rho_1^\# \pi_1^\#H = \rho_2^\# \pi_2^\#H = \pi^\#H.\]

        Furthermore, Hironaka's construction guarantees that
        \[
        \pi_{1}^{-1}(I(f_2)) \subset \pi_{1}^\# H \quad \text{and} \quad
        \pi_{2}^{-1}(I(f_1)) \subset \pi_{2}^\# H.
        \]
        So, we have
        \[\rho_1^* \pi_{2}^\#H = \pi^\#H + \sum_{i=1}^{r_2} a_{2i} F_{2i} \quad \text{and} \quad
        \rho_2^* \pi_{1}^\#H = \pi^\#H + \sum_{i=1}^{r_1} a_{1i} F_{1i}.\]
        Apply $\rho_l$ to $\phi_l^*H$ and get
        \[
        \widetilde{\phi}_1^*H
         = \rho_1^*\left( d_1 \cdot \pi_{1}^\#H \right) + \rho_1^*\left(\sum_{i=1}^{r_1} b_{1i} E_{1i} \right)
          =  d_1 \cdot \left(\pi^\#H +\sum_{j=1}^{r_2} a_{2j} F_{2j}\right)
            + \sum_{i=1}^{r_1} b_{1i} F_{1i}
        \]
        and
        \[
        \widetilde{\phi}_2^*H
         = \rho_2^*\left( d_2 \cdot \pi_{2}^\#H \right)+ \rho_2^*\left(\sum_{j=1}^{r_2} b_{2j} E_{2j} \right)
         =  d_2 \cdot \left(\pi^\#H +\sum_{i=1}^{r_1} a_{1i} F_{1i}\right)
            + \sum_{j=1}^{r_2} b_{2j} F_{2j}.
        \]
        Therefore,
        \begin{eqnarray*}
        \sum_{l=1}^2 \dfrac{1}{d_l}\widetilde{\phi}_l^*H  - \pi^*H
                          & = &  \sum_{l=1}^2 \left[ \left(\pi_V^\#H +\sum_{k\neq l} \sum_{j=1}^{r_k} a_{kj} F_{kj}\right)  + \dfrac{1}{d_l} \sum_{i=1}^{r_l} b_{li} F_{li} \right]\\
                          &   & - \left( \pi_V^\#H + \sum_{k=1}^2 \sum_{i=1}^{r_l} a_{li} F_{li} \right) \\
                          & = &  \pi_V^\#H + \sum_{l=1}^2 \left( \dfrac{1}{d_l} \sum_{i=1}^{r_l} b_{li} F_{li} \right)\\
                          & \succ &  \pi_V^\#H + \sum_{l=1}^2  \sum_{i=1}^{r_l} \dfrac{1}{r(\widetilde{\phi}_l)} a_{li} F_{li}
                          \quad \left( \because r(\widetilde{\phi}_l) \geq d_l \dfrac{a_{li}}{b_{li}} \right)\\
                          &\succ&  \min \left( \dfrac{1}{r(\widetilde{\phi}_l)} \right) \left(  \pi_V^\#H + \sum_{l=1}^2  \sum_{i=1}^{r_l} a_{li} F_{li} \right) \\
                          &=& \min \left( \dfrac{1}{r(\widetilde{\phi}_l)} \right) \pi^* H
        \end{eqnarray*}

         So, we get that the height function corresponding to an $\af^n$ effective divisor
         \[
         \sum_{l=1}^2 \dfrac{1}{d_l}\widetilde{\phi}_l^*H  - \left( 1 + \min \left( \dfrac{1}{r(\widetilde{\phi}_l)} \right) \right) \pi^*H
         \]
         is bounded below on $\pi^{-1}(\af^n)$ and hence get the desired result.
         \end{proof}

         \begin{cor}
        Let $f:\af^n \rightarrow \af^n$ be a regular polynomial automorphism defined over a number field $K$. Then, there is a constant $C$ such that
        \[
        \dfrac{1}{\deg f} h\bigl( f(P) \bigr) + \dfrac{1}{\deg f^{-1}} h\bigl( f^{-1}(P) \bigr)+ C
        > \left( 1+  \left( \dfrac{1}{\deg f \cdot \deg f^{-1}} \right) \right)h(P)
        \quad
        \text{for all }P \in \af^n(\overline{K}).\]
         \end{cor}
         \begin{proof} It is enough to show that $r(f)=r(f^{-1})=\deg f \cdot \deg f^{-1}$. \cite[Lemma 3.5]{Le} shows that
            \[{\pi_V}^*H = H_V + d'E_V + M_V\]
            \[{\phi_V}^*H = dH_V + E_V + I_V\]
        where $d$ is the degree of $\phi$, $d'$ is the degree of $\psi$ and $d'I - M$ is an effective divisor. Furthermore, since support of $d'I-M$ is not contained and doesn't contain $H_V$, it is actually $\af^n$-effective so that
        \[r(\phi) = d \times \max\left( \dfrac{1}{d}, d', \dfrac{M_i}{I_i}\right)\]
        where $M_i \leq d' I_i$. Since the support of $M_i, I_i$ does not contains $H$ so that $M_i \prec d' I_i$.
        Therefore,
        \[r(\phi)  = d \times d'\]
        \end{proof}

        \begin{ex} Let
                \[f(x,y) = (x^2+y,y).\]
               We can check that $r(f) = 1 \times \deg f = 2 $. Therefore, with any $g\in \Rat^2(H)$ which makes jointly regular pair with $f$ and $r(g)<2$, there is a constant $C$ such that
              \[\dfrac{1}{2} h\bigl( f(P) \bigr) + \dfrac{1}{d} h\bigl(g(P)\bigr) +C > \left (1+ \dfrac{1}{2}  \right) h(P)
              \quad \text{for all}~P \in \af^2(\overline{\mathbb{Q}}).\]
        \end{ex}

        \begin{rem}
            If a pair of rational maps is not jointly regular, then we may not have a similar inequality. For example, let
            \[
            f(x,y) = (x-y^d,y) \quad \text{and} \quad f^{-1}(x,y) = (x-y^d,y)
            \]
            where $d$ is a natural number larger than $2$. Then, for any $x$,
            \[
            \dfrac{1}{d}h\bigl( f(x,0) \bigr) + \dfrac{1}{d}h \bigl( f^{-1}(x,0) \bigr) = \dfrac{2}{d}h\bigl( (x,0) \bigr)
            \]
            so that it can't bound $h(x,0)+C$ above for any constant $C$.
        \end{rem}

    If $f_1\in \Rat^n(H)$ is a polynomial map such that $r(f_1) \geq \deg f_1$, then we do not get information from Theorem~A. However, if we can find a polynomial map $f_2 \in \Rat^n(H)$ such that $\{f_1,f_2\}$ is a jointly regular pair, then we can apply Theorem~B to get information for preperiodic points of a monoid generated by $\{f_1,f_2\}$.

    For each $m\geq 0$, let $W_m$ be the collection of ordered $m$-tuples chosen from
    $\{1, 2\}$,
    \[
    W_m = \bigl\{ (i_1, \cdots, i_m) ~|~ i_j \in \{1, 2\} \bigr\}
    \]
    and let
    \[
    W_* = \bigcup_{m\geq 0} W_m.
    \]
    Thus $W_*$ is the collection of words of $r$ symbols.

    For any $I=(i_1, \cdots, i_m) \in W_m$, let $f_I$ denote the composition of corresponding polynomial maps in $S$:
    \[
    f_I := f_{i_1} \circ \cdots \circ f_{i_m}.
    \]

    \begin{df}
        We denote \emph{the monoid of rational maps generated by $S=\{f_1, f_2\}$ under composition} by
        \[
        \Phi_S = \Phi := \{\phi = f_I ~|~ I \in W_*\}.
        \]
        Let $P\in \af^n$. \emph{The $\Phi$-orbit of $P$} is defined to be
        \[
        \Phi(P) = \{ \phi(P) ~|~ \phi \in \Phi \}.
        \]
        \emph{The set of (strongly) $\Phi$-preperiodic points }is the set
        \[
        \Pre(\Phi) = \{P \in \af^n ~|~ \Phi(P)~\text{is finite}\}.
        \]
    \end{df}

     \begin{thmE}
        Let $S=\{f_1, f_2 \}$ be a jointly regular pair of polyomial maps, let $f(f_l)$ be the $D$-ratio of $f_l$ and let $\Phi_S$ be the monoid of rational maps generated by $S$.
        Define
        \[
        \delta_S : = \left(\dfrac{1}{1+1/r} \right) \left(\dfrac{1}{\deg f_1} + \dfrac{1}{\deg f_2} \right)
        \]
        where $r = \displaystyle  \max_{l=1,2} \bigl(r(f_l)\bigr)$.

        If $\delta_S <1$, then
        \[
        \Pre(\Phi_S):= \bigcap_{f \in \Phi_S} \Pre(f) \subset \af^n_{\overline{K}}
        \]
        is a set of bounded height.
    \end{thmE}
    \begin{proof}
     By Theorem~B, we have a constant $C$ such that
    \begin{equation}\label{ineq1}\tag{G}
    0 \leq  \left(\dfrac{1}{1+\frac{1}{r}}\right) \sum_{l=1}^2\dfrac{1}{d_l} h\bigl( f_l(Q) \bigr) - h(Q) +C \quad \text{for all}~Q\in \af^n(\overline{K}).
    \end{equation}
    Note that if $r=\infty$, then $\left(\dfrac{1}{1+\frac{1}{r}}\right) = 1$, then it is done because of \cite[Theorem~4]{S4}. Thus, we may assume that $r$ is finite.

    We define a map $\mu:W_* \rightarrow \mathbb{Q}$ by the following rule:
    \[
    \mu_I = \mu_{(i_1, \cdots, i_m)} = \prod d_l^{p_{I,l}}
    \]
    where $p_{I,l} = - |\{t~|~ i_t = l\}|$. Then, by definition of $\delta_S$ and $\mu_I$, the following is true:
    \[
    \delta_S^m = \left[ \left(\dfrac{r}{r+1} \right)\sum_{l=1}^2 \dfrac{1}{d_l} \right]^m
    = \left(\dfrac{r}{r+1} \right)^m \sum_{I\in W_m} \dfrac{1}{\deg f_{i_1} \cdots \deg f_{i_m}}
    = \left(\dfrac{r}{r+1} \right)^m \sum_{I\in W_m} \mu_I.
    \]

    Let $P\in \af^n(\overline{\mathbb{Q}})$. Then, (\ref{ineq1}) holds for $f_I(P)$ for all $I\in W_m$:
    \begin{eqnarray*}
    0 \leq  \left(\dfrac{r}{r+1}\right) \sum_{l=1}^2 \dfrac{1}{d_l} h\bigl( f_l(f_I(P) ) \bigr) - h(f_I(P)) +C.
    \end{eqnarray*}
    Hence
    \begin{equation}\label{mainineq}\tag{H}
    0 \leq  \sum_{m=0}^M \sum_{I\in W_m}\mu_I \left(\dfrac{r}{r+1}\right)^m \left[ \sum_{l=1}^2 \dfrac{1}{d_l} h\bigl( f_l(f_I(P) ) \bigr) -
    \left(1+ \dfrac{1}{r}\right) h(f_I(P)) +C \right].
    \end{equation}

    The main difficulty of the inequality is to figure out the constant term. From the definition of $\delta_S$, we have
    \begin{eqnarray*}
    \sum_{m=0}^{M-1} \left(\dfrac{r}{r+1}\right)^m \sum_{I\in W_m} \mu_I = \sum_{m=1}^M \delta_S^m \leq \dfrac{1}{1-\delta_S}.
    \end{eqnarray*}
    Now, do the telescoping sum and most terms in (\ref{mainineq}) will be canceled:
    \begin{eqnarray*}\label{identity}
     &&\left(\sum_{m=0}^{M-1} \sum_{I \in W_m} \left(\dfrac{r}{r+1}\right)^m\mu_I \sum_{l=1}^2 \dfrac{1}{d_l} h\bigl( f_lf_I(P) \bigr)\right)
     - \left(\sum_{m=1}^{M} \sum_{I \in W_m} \left(\dfrac{r}{r+1}\right)^{m-1}\mu_I h\bigl( f_I(P) \bigr) \right)\\
     &&= \left(\sum_{m=0}^{M-1} \sum_{I \in W_m} \left(\dfrac{r}{r+1}\right)^m \sum_{l=1}^2 \dfrac{\mu_I}{d_l} h\bigl( f_lf_I(P) \bigr)\right)
     - \left(\sum_{m=0}^{M-1} \sum_{I \in W_m} \sum_{l=1}^2 \left(\dfrac{r}{r+1}\right)^m \dfrac{\mu_{I}}{d_l} h\bigl( f_lf_I(P) \bigr)\right) \\
     &&=0.
    \end{eqnarray*}
    Therefore, the remaining terms in (\ref{mainineq}) are the first term when $m=M$ and the last term when $m=0$. Thus, we get
    \begin{eqnarray*}
    0 &\leq&  \left[ \sum_{I\in W_M} \left(\dfrac{r}{r+1}\right)^M \mu_I \sum_{l=1}^K \dfrac{1}{d_l} h\bigl( f_l(f_I(P) ) \bigr)
    \right] - h(P) +  \sum_{I\in W_M} \left(\dfrac{r}{r+1}\right)^M \mu_I C \\
    &\leq&  \left[ \sum_{I\in W_M} \left(\dfrac{r}{r+1}\right)^M \mu_I \sum_{l=1}^2 \dfrac{1}{d_l} h\bigl( f_l(f_I(P) ) \bigr)
    \right] - h(P) + \dfrac{1}{1-\delta_S} C.
    \end{eqnarray*}
        Let $P$ be a $\Phi$-periodic point and define the height of the images of $P$ by the monoid $\Phi$ to be
   \[
   h(\Phi(P)) = \sup_{R\in \Phi(P)} h(R).
   \]
    Since
   \begin{equation*}
   \sum_{I\in W_M} \left(\dfrac{r}{r+1}\right)^M \mu_I \sum_{l=1}^2 \dfrac{1}{d_l} = \left(\dfrac{r}{r+1}\right)^M \sum_{I\in W_{M+1}} \mu_I
   = \left(1+\dfrac{1}{r}\right) \delta_S^{M+1},
   \end{equation*}
   and
   \[
   h(\Phi(P)) \geq h\bigl( g(P) \bigr) ~\text{for all}~g\in \Phi,
   \]
   we get
    \begin{eqnarray*}
    h(P) &\leq&  \left[ \sum_{I\in W_M} \left(\dfrac{r}{r+1}\right)^M \mu_I \sum_{l=1}^2 \dfrac{1}{d_l}\right]  h\bigl( \Phi(P) \bigr)+ \dfrac{1}{1-\delta_S} C\\
    &\leq&  \left(1+\dfrac{1}{r}\right) \delta_S^{M+1}  h\bigl( \Phi(P) \bigr)+ \dfrac{1}{1-\delta_S} C.
    \end{eqnarray*}
    By assumption, $\delta_S<1$ and $h \bigl(\Phi(P)\bigr)$ is finite, so letting $M \rightarrow \infty$ shows that
    $h(P)$ is bounded by a constant that depends only on $S$.
    \end{proof}
\section{Relation with Height expansion coefficient}

    In this section, we will discuss the relation between the $D$-ratio and the height expansion coefficient of dominant rational maps and find a lower bound of $r(f)$.

     Silverman defined the height expansion coefficient \cite[Definition 1]{S3} as follows:
   \begin{df}
    Let $\phi :W \dashrightarrow V$ be a rational map between quasiprojective varieties, all defined over $\overline{\mathbb{Q}}$. Fix height functions $h_V$ and $h_W$ on $V$ and $W$ respectively, corresponding to ample divisors. \emph{The height expansion coefficient of $\phi$} (relative to chosen height function $h_V$ and $h_W$) is the quantity
    \[
    \mu (\phi) = \sup_{\emptyset \neq U \subset W} \liminf_{\substack{P\in U \\ h(P) \rightarrow \infty }} \dfrac{h_V \bigl( \phi(P) \bigr) }{h_W(P)},
    \]
    where the sup is over all nonempty Zariski dense open subsets of $W$.
   \end{df}

       It seems that there is no direct relation between $\mu(f)$ and $r(f)$. So, define a new quantity which will be a bridge between them.

    \begin{df}
        Let $f\in \Rat^n(H)$ be a rational map defined over a number field $K$. Then, we define \emph{the height expansion coefficient of $f$ on $\af^n$}:
        \[
        c(f) := \liminf_{\substack{P\in \af^n \\ h(P) \rightarrow \infty} } \dfrac{h\bigl( f(P) \bigr)}{h(P)}.
        \]
    \end{df}
    Now, we can find a relation between $c(f)$ and $r(f)$.
    \begin{thm}
        \[c(f)^{-1} \leq  \dfrac{r(f)}{\deg f}.\]
    \end{thm}
    \begin{proof}  We may assume that  $P\in \af^n$ has sufficiently large $h(P)$.
    From Theorem~A, we have
    \[ \dfrac{r(f)}{\deg f}h\bigl( f(P) \bigr) > h(P) +C \quad
    \text{for all }P\in \af^n\]
    so that we may assume that $h\bigl( f(P) \bigr)$ is sufficiently large, too.
    Hence,
    \[ \dfrac{r(f)}{\deg f} > \dfrac{h(P) +C}{h\bigl( f(P) \bigr)} \] for $P\in \af^n$ with sufficiently large $h(P)$.
    Therefore,
    \[ \dfrac{r(f)}{\deg f}  \geq \limsup_{h(P) \rightarrow \infty} \dfrac{h(P) +C}{h\bigl( f(P) \bigr)} = \dfrac{1}{c(f)}.\]
    \end{proof}

    \begin{cor}
        Let $f\in \Rat^n(H)$ be a rational map defined over a number field $K$. If there is a curve $C$ on $\af^n$ whose image under $f$ is a point, then $r(f) = \infty$.
    \end{cor}
    \begin{proof} It is enough to show $c(f)^{-1}=\infty$. Let $P_n$ be a sequence on the given curve $C$ whose height goes to infinity. Then,
    \[
    c(f)^{-1} = \limsup_{h(P) \rightarrow \infty} \dfrac{h(P)}{h\bigl( f(P) \bigr)}
    \geq \lim_{n \rightarrow \infty} \dfrac{h(P_n)}{h\bigl( f(P_n) \bigr)}
    = \lim_{n \rightarrow \infty} \dfrac{h(P_n)}{h(Q)} = \infty\]
    \end{proof}

    Furthermore, the relation between $c(f)$ and $\mu(f)$ is clear and hence we can build the following inequality.
    \begin{prop}\label{mu(f) r(f)}
    Let $f\in \Rat^n(H)$ be a rational map defined over a number field $K$. Then
        \[
        \mu(f) \geq c(f) \geq \dfrac{\deg f}{r(f)}.
        \]
    \end{prop}
    \begin{proof}
        It is clear since
        \[
        c(f) = \liminf_{\substack{P\in \af^n \\ h(P) \rightarrow \infty} } \dfrac{h\bigl( f(P) \bigr)}{h(P)}.
        \]
        is the case when $U = \af^n$.
    \end{proof}
    With Proposition~\ref{mu(f) r(f)}, we can find an easier proof of Proposition~\ref{D-ratio prop}~(1).
    \begin{cor}
    Let $f\in \Rat^n(H)$ be a rational map defined over a number field $K$. Then
        \[
        r(f) \geq 1.
        \]
    \end{cor}
    \begin{proof}
    Remind that we have a lower bound of $h(P)$: by \cite[Theorem~B.2.5]{SH}, there exists a constant $C$ such that
    \[
    \deg f \cdot h(P) +C  > h\bigl( f(P) \bigr) \quad \text{for all}~P \in \pp^n(\overline{K}) \setminus I(f).
    \]
    So, for any open set $U$ of $\pp^n$, we have
    \[
    \liminf_{\substack{P\in U \\ h(P) \rightarrow \infty }} \dfrac{h \bigl( f(P) \bigr) }{h(P)}
    \leq \liminf_{\substack{P\in U \\ h(P) \rightarrow \infty }} \dfrac{\deg f \cdot h(P) +C  }{h(P)} = \deg f.
    \]
    and hence
    \[
    \mu (f) = \sup_U \liminf_{\substack{P\in U \\ h(P) \rightarrow \infty }} \dfrac{h \bigl( f(P) \bigr) }{h(P)}\leq  \deg f.
    \]

    Therefore, by Proposition~\ref{mu(f) r(f)}, we get
    \[
    \deg f \geq \mu(f) \geq \dfrac{\deg f}{r(f)}.
    \]
    \end{proof}

    \begin{ex}
        Let $f [x,y,z]=[x^{km}, y^{(k-1)m} z^{m},z^{km}]$. \\

        First, since $f \equiv [1,0,0]$ on $H \setminus [0,1,0]$,
        \[
        \liminf_{\substack{h(P) \rightarrow \infty \\ P\in U' }}\dfrac{h\bigl(f(P)\bigr)}{h(P)}=0
        \]
        for all open set $U' \not\subset \af^2$. Thus,
        \[
        \mu(f) = \sup_{U: open set \neq \emptyset} \liminf_{\substack{h(P) \rightarrow \infty \\ P\in U }}\dfrac{h(f(P))}{h(P)}
        = \sup_{U \subset \af^2, open } \liminf_{\substack{h(P) \rightarrow \infty \\ P\in U }}\dfrac{h(f(P))}{h(P)}.
        \]
        So, we may assume that $U \subset \af^2$.

        Let $T_\alpha = \{[x,y,z] ~|~ x = \alpha z \}$.
        Then, $\displaystyle \bigcup_{\alpha\in \mathbb{\mu}}$ is Zariski dense in $U$: there is an $\alpha$ such that
        \[
        U  \cap T_\alpha \neq \emptyset.
        \]
        Because $U  \cap T_\alpha$ is an open set of $T_\alpha$, there is a sequence $P_M= [\alpha ,y_M,1] \in U \cap T_\alpha$ such that $\displaystyle\lim_{M \rightarrow \infty} h(P_M) = \infty$. From the triangle inequality and the definition of the height, we get
        \[
        h(y_M) \leq h(P_M) \leq h(y_M) + h(\alpha) = h(y_M)
        \]
        and
        \[
        h\left(y_M^{(k-1)m} \right) \leq h \bigl( f(P_M) \bigr) = \left[\alpha^d, y_M^{(k-1)m},1\right] \leq  h\left(y_M^{(k-1)}m\right).
        \]
        Moreover, if $\displaystyle\lim_{M \rightarrow \infty} h(P_M) = \infty$, then
        \[
        \lim_{M\rightarrow \infty} \dfrac{h \bigl( f(P_M) \bigr)}{h(P_M)}
        = \lim_{M\rightarrow \infty} \dfrac{(k-1)m \cdot h(y_M) }{ h(y_M)} = (k-1)m
        \]
        and hence
        \begin{eqnarray*}
        \liminf_{\substack{h(P) \rightarrow \infty \\ P\in U }} \dfrac{h\bigl(f(P)\bigr)}{h(P)} \leq (k-1)m.
        \end{eqnarray*}

        Next, figure out lower bound with $r(f)$: we have $(V, \pi_V)$ a resolution of indeterminacy of $f$ by $k$ successive blowups and get a resolved morphism $\phi_V$. Then, pull-backs of $H$ by $\pi_V$ and $\phi_V$ are calculated as follows;
        \[
        \pi_V^*H = H_V + \sum_{i=1}^k iE_i \quad \text{and} \quad \phi_V^*H = km H_V + \sum_{i=1}^k im(k-1)E_i.
        \]
        Therefore,
        \[
        r(f) = d \cdot \max_i \left( \dfrac{1}{d}, \dfrac{i}{im(k-1)} \right) = km\dfrac{1}{m(k-1)} = \dfrac{k}{k-1}.
        \]
        and hence
        \[
        (k-1)m = \mu(f) \geq c(f) \geq \dfrac{\deg f}{r(f)} = d\dfrac{k-1}{k} = (k-1)m.
        \]
    \end{ex}


\begin{thebibliography}{ABCDEFG}

        \bibitem{CS}
            Call, Gregory S.; Silverman, Joseph H. Canonical heights on varieties with morphisms. Compositio Math. 89 (1993), no. 2, 163-205.

        \bibitem{Cu}
             Cutkosky, Steven Dale. {\em Resolution of singularities},
             Graduate Studies in Mathematics, Vol 63, American Mathematics Society, 2004.

        \bibitem{De}
             Demailly, J.-P. , {\em Multiplier ideal sheaves and analytic methods in algebraic geometry, in: School on Vanishing Theorems
        and Effective Results in Algebraic Geometry}, Trieste 2000, ICTP Lecture Notes Vol. 6 (Abdus Salam Int. Cent. Theoret.
        Phys., Trieste, 2001), pp. 1-148.


        \bibitem{D}
            Denis, L.,  {\em Points p\'{e}riodiques des automorphismes affines}, J. Reine Angew. Math. 467, 157-167, 1995

        \bibitem{GJ}
             Fisher, G., Piontkowski, J. {\em Ruled Variety, an introduction to algebraic differential geometry},
             Advanced Lectures in Mathematics. Friedr. Vieweg \& Sohn, Braunschweig, 2001.

        \bibitem{F}
             Fulton, W., {\em Intersection theory },  Second edition, Springer-Verlag, Berlin, 1998.

        \bibitem{H}
            Hartshorne, S., {\em Algebraic geometry}, Springer, 1977.

        \bibitem{Hi}
            Hironaka, H., {\em Resolution of singularities of an algebraic variety over a field of characteristic zero. I}, Ann. of Math. (2) 79(1964), 109-203.

        \bibitem{HPV}
            Hubbard, J., Papadopol, P., Veselov, V. {\em A compactification of H\'{e}non mappings in $\mathbb{C}^2$ as dynamical systems}.
            {\em Acta Math.}, 184:203-270, 2000.

        \bibitem{K}
             Kawaguchi, S., {\em Canonical height functions for affine plane automorphisms,} {\em Math. Ann.}, 335, no. 2(2006), 285-310.

        \bibitem{K2}
             Kawaguchi, S., {\em Local and global canonical height functions for affine space regular automorphisms}, preprint, Arxiv:0909.3573, 2009.

        \bibitem{La}
            Lang, S., {\em Fundamentals of diophantine geometry}, Berlin Heidelberg New York: Springer 1983.

        \bibitem{Laz}
        Lazarsfeld, R. , {\em Positivity in Algebraic Geometry I}, Ergebnisse der Mathematik und ihrer Grenzgebiete 3. Bd. 48 (Springer, New York, 2004)

        \bibitem{Le}
            Lee, C., {\em An upper bound for height for regular affine automorphisms on $\af^n$}, preprint, ArXiv:0909.3107, 2009.


       \bibitem{M}
            Marcello, S., {\em Sur les propietes arithmetiques des iteres d'automorphismes reguliers}, C. R. Acad. Sci. Paris Ser. I Math. 331, 11-16., 2000.

       \bibitem{MS}
            Morton, P, Silverman, J.H, {\em Rational periodic points of rational functions},
            Internat. Math. Res. Notices, (2) (1994), 97-110.
        \bibitem{N}
            Northcott, D. G., {\em Periodic points on an algebraic variety}, Ann. of Math. (2), 51(1950), 167-177.

        \bibitem{Sh}
            Shafarevich, I. {\em Basic Algebraic Geometry}, Springer, 1994.

        \bibitem{Si}
            Sibony, N. {\em Dynamique des applications rationnelles de $\pp^k$}, {\em Panor. Syntheses}, 8(1999), 97-185.

        \bibitem{S1}
            Silverman, J. H. {\em Geometric and arithmetic properties of the Henon map}, {\em Math. Z.}, 215, no. 2(1994), 237-250.

        \bibitem{S4}
            Silverman, J. H. {\em Height bounds and preperiodic points for families of jointly regular affine maps.}, Pure Appl. Math. Q., 2, no. 1, part 1(2006), 135-145.

        \bibitem{S3}
            Silverman, J. H. {\em Height estimate for equidimensional dominant rational maps}, preprint, ArXiv:0908.3835, 2009.

        \bibitem{SH}
            Silverman, J. H., Hindry, M. {\em Diophantine Geometry, An introduction}, Springer, 2000.

        \bibitem{S2}
            Silverman, J. H. {\em The arithmetic of Dynamical system}, Springer, 2007.

    \end{thebibliography}
\end{document}